\documentstyle[12pt]{article}
\renewcommand {\c}  {\'{c}}
\newcommand {\cc} {\v{c}}
\newcommand   {\s}  {\v{s}}

\newcommand{\somk}{\sum_{k=1}^{\infty}}
\begin{document}

%
%**************************TITLE PAGE***************************************
%

\begin{titlepage}
\baselineskip=24pt
\title{UNIFIED VIEW OF MULTIMODE ALGEBRAS WITH FOCK-LIKE REPRESENTATIONS}
\author{Stjepan Meljanac and Marijan Milekovi\c $^\dagger$ }
\maketitle
\bigskip
\begin{center}
{\it Rudjer Bo\s kovi\c \ Institute\\Bijeni\cc ka c.54, 41001 Zagreb,
Croatia}\\
{\it $^\dagger$ Prirodoslovno-Matemati\cc ki Fakultet,Zavod za teorijsku
fiziku,\\
Bijeni\cc ka c.32, 41000 Zagreb, Croatia}\\
\bigskip
Short title: {\it Unified view of multimode...}\\
( to be published in Int.J.Mod.Phys. A )

\end{center}
\end{titlepage}
%%%%%%%%%%%%%%%%%%%%%%%%%%%%%%%%%%%%%%%%%%%%%%%%%%%%%%%%%%%%%%%%%%%%%%%%%%%
%%%%%%%%%%%%%%%%%%%%%%%%%%%%%%%ABSTRACT%%%%%%%%%%%%%%%%%%%%%%%%%%%%%%%%%%%%

\begin{center}
{\large \bf Abstract}
\end{center}
\vspace{0.5cm}
\baselineskip=24pt
A unified view of general multimode oscillator algebras
with Fock-like
representations is presented.It extends a previous
analysis of the single-mode oscillator
algebras.The expansion of the $a_ia_j^{\dagger}$
operators is extended
 to include all normally ordered terms in
  creation and annihilation operators and we
  analyze their action on Fock-like
  states.We restrict ourselves to the algebras
  compatible with number operators.
  The connection between these algebras and generalized
  statistics is analyzed.We demonstrate our approach by
  considering the algebras
   obtainable
  from the generalized Jordan-Wigner transformation,
  the para-Bose and
  para-Fermi algebras,
  the Govorkov "paraquantization" algebra and generalized
  quon algebra.

%%%%%%%%%%%%%%%%%%%%%%%%%%%%%%%%%%%%%%%%%%%%%%%%%%%%%%%%%%%%%%%%%%%%%%%%%%%
%%%%%%%%%%%%%%%%%%%%%%%%%%%%%INTRODUCTION%%%%%%%%%%%%%%%%%%%%%%%%%%%%%%%%%%
\newpage
\section{Introduction}

Recently,much attention has been devoted to the study of quantum groups and
algebras$^1$,
noncommutative spaces and geometries$^2$,generalized notions of symmetries as
well as
to their diverse applications in physics$^3$. These
approaches are in close relationship with study of deformed
oscillators,algebras and their
 Fock representations. Single-mode oscillator algebras were studied by  a
 number of authors$^4$.
 A  unified view of single-mode oscillator algebras was   proposed by
 Bonatsos and Daskaloyannis $^5$ and
 Meljanac et al.$^6$. Multimode oscillator algebras are much more complicated
 and only partial results exist in the literature$^7$.Particularly,the R-matrix
 approach to   multimode algebras was followed in Refs.$(8)$.Some of
 multimode algebras,but not all, can be obtained from an
  ordinary Bose
 algebra with equal number of oscillators$^{9}$.\\
 On the other hand, there is an old additional physical motivation to study
multimode oscillator
 algebras,connected with the problem of generalized statistics,different from
 Bose and Fermi statistics$^{10}$.The first consistent example of it was
para-Bose and
 para-Fermi statistics$^{11,12}$ and, recently, a new paraquantization$^{13}$
satisfying
 trilinear commutation relations between annihilation and creation operators.
 These types of
 statistics are characterized by a discrete parameter called the order
 of parastatistics ${\it p} \in {\bf N}$,interpolating between Bose and Fermi
statistics.\\
 Recently,a new interpolation,namely  infinite quon statistics characterized by
a
 continuous parameter, has been proposed and analyzed$^{14,15}$.Its
multiparameter
 extension was studied in Refs.(16).Parastatistics can be applied in spaces
with
  an arbitrary number of dimensions and those statistics with a continuous
parameter
 could be relevant to lower dimensions.Specially,in the $(2+1)$ dimensional
  space,another type is also  possible,
 the so-called anyonic (fractional) statistics
 $^{17}$.It is characterized by the continuous statistical parameter $\lambda$.
 A different type of the fractional statistics,generalizing the Pauli
 exclusion principle,is also proposed$^{18}$.\\
  Our motivation of the present work is twofold:\\
  (i) to try to extend the quantization defined by $a_ia^{\dagger}_j-R_{ij,kl}
  a^{\dagger}_ka_l=\delta_{ij}$,$i,j \in I$ to include all normally ordered
terms in
  creation and annihilation operators;\\
  (ii) to connect these algebras with the notion of generalized statistics.\\
In the physical world particles decay to other particles,interact and transform
 themselves,and, as a consequence,the number operators in these processes are
not
 conserved.They change respecting some selection rules,characteristic of the
 type of interaction.However,in order to measure different sorts of particles,
 we always assume that any sort of physical particles is countable.Thus ,the
existence
 of the number operator is of utmost importance.
 Our aim is to define and describe the most general multimode oscillator
algebras
  possessing Fock-like representations and compatible number
operators.Generally,
  one considers a given algebra and examines all its representations.Here we
are going
  backwards,i.e. we start with the most general Fock-like space and look for
the relations
between annihilation and creation operators $\{ a_i$,$\bar{a}_i|i \in
I\}$,leading
to the positive definite scalar product and to the non-negative norms of all
Fock-like states. \\
 The plan of the paper is as follows.In Sec. 2 we briefly review the main
  aspect of a
single-mode oscillator algebra.We characterize the algebra by a set
of quantities $\{\varphi(n)\}$,$\{c_n\}$ and $\{d_n\}$ which are
 generalized to a set $\{\Phi\}$,$\{C\}$ and $\{D\}$ for
multimode oscillator algebra with Fock-like representations in  Sec. 3.
We also discuss the connection with generalized statistics.
In Sec. 4-7 we give some examples of our approach to  multimode algebras, i.e.
we discuss algebras obtainable from  generalized Jordan-Wigner transformations
(Sec.4),para-algebras (Sec.5),new paraquantization of Govorkov (Sec.6)
and generalized quon algebra (Sec.7).

%%%%%%%%%%%%%%%%%%%%%%%%%%%%%%%%%%%%%%%%%%%%%%%%%%%%%%%%%%%%%%%%

\section{Single-Mode Oscillator Algebra with Fock-like Representations}
\setcounter{equation}{0}

In this section we briefly review  single-mode oscillator algebras.For
more details we refer to Ref.(6).\\
Let us consider a pair of operators $\bar{a}$, $a$ (not necessarily Hermitian
conjugate to each other) with the number operator N. The most general
commutation relation linear in the $\bar{a} a$ and $a \bar{a}$ operators is
\begin{equation}
  a \bar{a} - F(N) \bar{a} a = G(N),
\end{equation}
where $F(N)$ and $G(N)$ are arbitrary complex functions. The number
operator satisfies
\begin{equation}
\begin{array}{c}
[N,a] = -a,\\[4mm]
[N,\bar{a}] = \bar{a},\\[4mm]
[N, \bar{a} a] = [N,a \bar{a}] = 0 .
\end{array}
\end{equation}
Hence, we can write
\begin{equation}
\begin{array}{c}
\bar{a} a = \varphi (N),\\[4mm]
a \bar{a} =  \varphi (N+1),
\end{array}
\end{equation}
where $\varphi (N)$ is, in general, a complex function satisfying
the recurrence relation
\begin{equation}
\varphi (N+1) - F(N) \varphi (N) = G(N).
\end{equation}
If $\varphi (N)$ is the bijective mapping, then
\begin{equation}
N = \varphi^{-1}(\bar{a} a) = \varphi^{-1}(a \bar{a}) - 1.
\end{equation}

Let us denote  the Hermitian conjugate of the operator
$a$ by  $a^{\dagger}$ . Then it follows that
\begin{equation}
\begin{array}{c}
[N,a^{\dagger}] = a^{\dagger},\\[4mm]
\bar{a}  = c(N) a^{\dagger},
\end{array}
\end{equation}
where $c(N)$ is a complex function of $N$. It is convenient to choose
$c(N)$ to be a "phase" operator, $|c(N)| = 1$. Then we have
\begin{equation}
\begin{array}{c}
a^{\dagger} a = |\varphi (N)|,\\[4mm]
a a^{\dagger} = |\varphi (N+1)|,\\[4mm]
a a^{\dagger} - a^{\dagger} a  =  |\varphi (N+1)| - |\varphi (N)| = G_{1}(N),
\\[4mm]
c(N) = e^{i \arg \varphi (N)} = \frac{\varphi (N)}{|\varphi (N)| }.
\end{array}
\end{equation}
If $\varphi (N) > 0$, then $\arg \varphi (N) = 0$ and $c(N) = 1$.\\
Let us further assume that $|0>$ is a vacuum:
 \begin{equation}
\begin{array}{c}
a |0> = 0,\\[4mm]
N |0> = 0 , \quad    \varphi (0) = 0,\\[4mm]
<0|0> = 1 .
\end{array}
\end{equation}
One can always normalize the operators $a$ and $\bar{a}$ such that
$|\varphi (1)| = 1$; then\\ $|G(0)| = 1$. The function $\varphi (N)$ is
determined by the recurrence relation (2.4) and is given by
\begin{equation}
\varphi (n) = [F(n-1)]! \sum_{j=0}^{n-1} \frac{G(j)}{[F(j)]!},
\end{equation}
where
\begin{equation}
\begin{array}{c}
[F(j)]! = F(j)F(j-1)...F(1),\\[4mm]
[F(0)]! = 1 .
\end{array}
\end{equation}
The excited states with unit norms are
$$
|n> = \frac{(a^{\dagger})^{n}}{\sqrt {[|\varphi (n)|]!}} |0> =
\frac{(c^{-\frac{1}{2}} \bar{a})^{n}}{\sqrt {[\varphi (n)]!}} |0>,\\[4mm]
$$
\begin{equation}
\begin{array}{c}
<n|m> = \delta_{mn},   \quad   n,m = 0,1,2....,\\[4mm]
<n-1|a|n> = <n|a^{\dagger}|n-1> = \sqrt {|\varphi (n)|} .
\end{array}
\end{equation}

The function $|\varphi (n)|$ uniquely determines the type of deformed
oscillator algebra and vice versa. If $\varphi _{1}$ $\neq$ $\varphi _{2}$
but
 $|\varphi _{1}| = |\varphi _{2}|$, the corresponding algebras are isomorphic.
There is a family of  functions $(F,G)$ leading to the same algebra, with
identical
 functions $\varphi (N)$. Therefore we can fix the "gauge", for example:
\begin{equation}
\begin{array}{c}
(a)  \quad       F(N) = 1 ,    \quad    a \bar{a} - \bar{a} a =
G_{1}(N),\\[4mm]
(b)   \quad     G(N) = 1,     \quad   a \bar{a} - F_{1}(N)\bar{a} a = 1,
\\[4mm]
(c)    \quad     F(N) = q,      \quad   a \bar{a} - q \bar{a} a = G_{q}(N) .
\end{array}
\end{equation}
Several examples are given in Ref.(6).\\
The unitary irreducible representations of the single-mode oscillator algebra
can be classified according to the existence (or nonexistence)
of a vacuum state $|0>$
 as Fock-like,non-Fock-like and degenerate.
If $\varphi (n) \neq 0$ for $\forall n$ $ \in$ ${\bf N}$, then there is an
infinite
set ("tower") of states and ,owing to the well-defined mapping to the Bose
algebra
\begin{equation}
a = b \sqrt{\frac{|\varphi (N)|}{N}},  \quad    a^{\dagger} = \sqrt{\frac
{|\varphi (N)|}
{N}} b^{\dagger},
\end{equation}
the Fock space of the deformed algebra is identical to the Fock space of the
Bose oscillator.
 However, if  $\varphi (n_{0}) = 0$ for some $n_{0}$,
 then the state $(a^{\dagger})^{n_{0}} |0>$ has zero norm and, consistently,
we can put $|n_{0}> \equiv  0$. The corresponding representation is
finite-dimensional and the representation matrices are of the
$n_{0} \times n_{0}$ type (2.11) (degenerate representation).

Now we present another approach to the single-mode oscillator algebras which
will be generalized to the multimode case in the next section.\\
One can write Eq.(2.3) in an alternative form
\begin{equation}
\varphi(N+1)=a\bar{a}=1+\sum_nc_n(\bar{a})^n(a)^n
\end{equation}
where the coefficients $c_n$ are defined recursively by ${\varphi (n)}$ as

\begin{equation}
c_n=\frac{\varphi (n+1)-1-\sum_{k=1}^{n-1}c_k\varphi(n)\cdots \varphi(n+1-k)}
{[\varphi(n)]!},\qquad \forall \varphi(n) \neq 0.
\end{equation}

Knowing $c_k$ we can recursively obtain $\varphi(n)$
\begin{equation}
\varphi(n+1)=\sum_{k=0}^n c_k \frac{[\varphi(n)]!}{[\varphi(n-k)]!},\qquad
c_0=1.
\end{equation}
Equivalently,one can write the expression for the number operator N,Eq.(2.5),as
\begin{equation}
N= \bar{a}a+\sum_{n\geq 2}d_n \, (\bar{a})^n(a)^n,,
\end{equation}
where
\begin{equation}
d_n=\frac{n-\sum_{k=1}^{n-1} d_k \, \varphi(n)\cdots \varphi(n+1-k)}
{[\varphi(n)]!}.
\end{equation}
$$d_1=1 ,\qquad \forall \varphi(n) \neq 0 .$$
Note that $d_n=0$ for $\varphi(n)=0$,$\varphi(n-1)\neq 0$.Knowing $d_n$, we can
obtain
$\varphi(n)$ from the same recurrent relation (2.18).\\
In a physical application it is important to know the vacuum matrix elements
\begin{equation}
A_{m,n}=<0|(a)^n(\bar{a})^n|0>=[\varphi(n)]!\delta_{mn}.
\end{equation}
Hence,the sets of quantities $\{\varphi(n)\}$,$\{c_n\}$ and $\{d_n\}$ represent
a
unique description of a given (deformed) single-mode oscillator algebra.
Together with the generalization of the matrix ${\bf A}$,Eq.(2.19),they will be
particularly
useful for the analysis of   multimode (deformed) oscillator algebras.

%%%%%%%%%%%%%%%%%%%%%%%%%%%%%%%%%%%%%%%%%%%%%%%%%%%%%%%%%%%%%%%%%

\section{Multimode Oscillator Algebra with Fock-like Representations}
\setcounter{equation}{0}
In this section we define and describe the most general multi-mode oscillator
algebras   possessing Fock-like representations.We
start with the most general Fock-like space and look for the relations
between annihilation and creation operators $\{ a_i$,$\bar{a}_i|i \in
I\}$,leading
to the positive definite scalar product and to the non-negative norms of all
Fock-like states (vectors).Furthermore,we  restrict ourselves to the algebras
possessing  number operators for all $\forall i \in I $.\\
Let us first construct the most general Fock-like space.Let there  be  a vacuum
state $|0\rangle$ and let us assume that there exist independent annihilation
and creation operators $\{ a_i$,$\bar{a}_i|i \in I\}$,satisfying the vacuum
condition
\begin{equation}
a_i|0\rangle  = 0, \qquad \forall i \in I .
\end{equation}
The set $\{I\}$ can be finite or infinite,discrete or continuous and with some
(partial) ordering.It may  have some additional mathematical and physical
structure.
\\Now,we build  a Fock-like space starting from the vacuum state $|0\rangle$.
The space $ {\cal{H}}_0=\{\lambda |0\rangle |\lambda \in {\bf C}\}$ is
one-dimensional.
We assume that the states $|0\rangle$ and $\bar{a}_i|0\rangle$ are linearly
independent states and they span the linear space ${\cal{H}}_1=
\{\sum_{i}\lambda_i \,\bar{a}_{i}|0\rangle|i\in I\}$  and ${\it dim}$
${\cal{H}}_{1} = |I|$.
Similarly,the states $\bar{a}_{i_1}\bar{a}_{i_2}|0\rangle$ span
 the linear space $${\cal{H}}_2=\{\sum_{i_1,i_2} \lambda_{i_1i_2}\bar{a}_{i_1}
\bar{a}_{i_2}|0\rangle | i_1,i_2 \in I \}$$ with the dimension $dim {\cal{H}}_2
\leq |I|^2$.
We assume that any state belonging to ${\cal{H}}_2$ is linearly independent of
the
states in ${\cal{H}}_0 \bigoplus {\cal{H}}_1$.Generally,
\begin{equation}
\begin{array}{c}
{\cal{H}}_n=\{\sum_{i_1\cdots i_n }\lambda_{i_1\cdots i_n }\, \bar{a}_{i_1}
\cdots \bar{a}_{i_n}
 |0\rangle \; |\;
\lambda_{i_1 \cdots i_n} \in {\bf C }\},\\[4mm]

{\cal{H}}_n\cap {\cal{H}}_m = \{0 \} ,\qquad n \neq m\\[4mm]

dim {\cal{H}}_n \leq |I|^n .
\end{array}
\end{equation}
Hence,the full Fock-like space ${\cal{H}}$ is
\begin{equation}
{\cal{H}} = {\cal{H}}_0 \bigoplus {\cal{H}}_1 \bigoplus {\cal{H}}_2 \bigoplus
\cdots
= \bigoplus_n {\cal{H}}_n,\qquad n \in {\bf N_0}.
\end{equation}
Furthermore,we assume that the annihilation operators $a_i,i \in I$ act on the
space
${\cal{H}}_n$ in such a way that
\begin{equation}
a_{i_1}\cdots a_{i_k}\bar{a}_{j_1} \cdots \bar{a}_{j_n}|0\rangle \in
{\cal{H}}_{n-k},
\qquad n \geq k,
\end{equation}
and this state is equal to zero if $n < k$.Then,the total number operator N can
be
defined as
\begin{equation}
N\, v_n=n\, v_n,\qquad v_n \in {\cal{H}}_n,\; n\in {\bf N_0}.
\end{equation}
${\cal{H}}_n$ is an invariant subspace with respect to N.We call the states
belonging to
 ${\cal{H}}_n$ ${\it n-particle}$ states.The number of independent n-particle
states
  is equal to ${\it dim}{\cal{H}}_n$ and is connected with the problem of
generalized
  statistics.We discuss this point at  the end of this section.\\
  We also  assume that there exist the number operators $N_i$,$i\in I$,for
 every species ${\it i}$.Namely,
 \begin{equation}
 \begin{array}{c}
 N_i (\bar{a}_{i_1}\cdots \bar{a}_{i_n})|0\rangle=\sum_{k=1}^{n} \delta_{ii_k}
 (\bar{a}_{i_1}\cdots \bar{a}_{i_n})|0\rangle,
 \qquad i_1 \cdots i_n \in I,\quad n\in{\bf N}\\[4mm]

 N_i|0\rangle=0,\qquad \forall i \in I ,\\[4mm]

 [N_i,\bar{a}_j]=\delta_{ij} \bar{a}_i, \qquad [N_i,a_j]=-\delta_{ij}
a_i,\\[4mm]

N=\sum_{i\in I}N_i,\qquad [N_i,N_j]=0,\qquad \forall i,j \in I.
\end{array}
\end{equation}
Hence the subspace $${\cal{H}}_{i_1 \cdots i_n}=\{\sum_{\pi \in S_n}
\lambda_{\pi \cdot (i_1 \cdots i_n)} \pi \cdot (\bar{a}_{i_1} \cdots
\bar{a}_{i_n})|0\rangle \;|
\; i_1 \leq i_2 \leq \cdots \leq i_n \in I\}
 \subset {\cal{H}}$$ is invariant  with  respect to N and ${\cal{H}}_{i_1
\cdots i_n}
 \cap {\cal{H}}_{j_1 \cdots j_n}= \{ 0 \}$ if $i_k \neq j_k$.\\
 Let us also remark that if we perform a mapping from  the operators
  $\{\bar{a}_i,a_i| i \in I\}$ to operators $\{\bar{a}_j',a_{j}'| j \in I'\}$
  \begin{equation}
  a_{j}'=f_j(a_i) \qquad \bar{a}_j'=\bar{f_j}(\bar{a_i})
  \end{equation}
  and repeat the same construction of the Fock-like space ${\cal{H}}'$ under
  the vacuum state $|0\rangle '\equiv |0\rangle $,then,generally,the number
operators $N_{j}'$
  will not exist.

 Now,let us define the dual Fock-like space ${\cal{H}}^{(d)}$,based on the dual
  vacuum-state vector
 $\langle 0|$,with the condition $\langle 0|0\rangle=1$.The construction of
 ${\cal{H}}^{(d)}$
 is the same as that for ${\cal{H}}$.Note that $\langle 0|\bar{a}=0$,$ \forall
i \in I$.
 The dual ${\it n-particle}$ states are
 \begin{equation}
 \langle 0 |\, (a_{i_1} \cdots a_{i_n})\; \in {\cal{H}}^{(d)}_{i_1 \cdots i_n}
 \end{equation}
   We also have, as in Eq.(3.4),
 \begin{equation}
 \begin{array}{c}
 \langle 0 | \, (a_{j_1} \cdots a_{j_n}\bar{a}_{i_1} \cdots \bar{a}_{i_k})\;
\in
  {\cal{H}}^{(d)}_{n-k} \qquad  n \geq k\\[4mm]
\langle 0 |  \, (a_{j_1} \cdots a_{j_n}\bar{a}_{i_1} \cdots \bar{a}_{i_k})=0
 \qquad  n<k .
 \end{array}
 \end{equation}
 From the existence of the number operators $N_i$,Eq.(3.6),it follows that
\begin{equation}
 {\cal{H}}^{(d)}_{i_1 \cdots i_n} \cap {\cal{H}}^{(d)}_{j_1 \cdots j_n}= \{ 0
\},
 \qquad  i_{k} \neq j_{k}.
\end{equation}
 Note that ,so far, we have assumed no precise relation between the
 $a_i$ and $\bar{a}_i$ operators.

 Our next step is to construct the scalar product and define the norm of any
vector
 in the Fock-like space.Therefore,we   construct a sesquilinear form
  $<*,* > : {\cal{H}}^{(d)}\times {\cal{H}}\rightarrow {\bf C}$.Namely,we need
  to calculate the matrix ${\bf A}$ with the matrix elements
  \begin{equation}
  A_{i_1\cdots i_n;j_1 \cdots j_m}=\langle 0 |a_{i_1} \cdots
a_{i_n}\bar{a}_{j_1}
   \cdots \bar{a}_{j_m}| 0 \rangle.
   \end{equation}
   From Eq.(3.6),i.e. the existence of the number operators $N_i$,it follows
that
   \begin{equation}
A_{i_1\cdots i_n;j_1 \cdots j_m}=0 \qquad  (i)\neq (j),
\end{equation}
i.e.it vanishes unless $n=m$ and the indices $(i_1\cdots i_n)$ and
 $(j_1 \cdots j_m)$ are equal up to permutation.\\
We can always choose
\begin{equation}
\langle 0 |0 \rangle=1,\qquad \langle 0 |a_i\bar{a}_j |0 \rangle=\delta_{ij}
\end{equation}
{}From the physical point of view,we need to know how to calculate any monomial
$P(a_i,\bar{a}_j)$ in $a_i,\bar{a}_j$ between vacuum states.Hence,in order to
calculate all such matrix elements, it is sufficient to know the action of all
annihilation operators $a_i$ on all Fock-like states
\begin{equation}
a_i\bar{a}_{i_1}\cdots \bar{a}_{i_n}|0\rangle \in {\cal{H}}_{ i_1\cdots
\not  i\cdots i_n} \subset
{\cal{H}}_{n-1},
\end{equation}
where the slash denotes the omission of the corresponding index (state) from
the
set
$i_1\cdots i_n$.
Hence,these relations are compatible with the number operators
$N_i$,Eq.(3.6),if
and only if
$$
a_i\bar{a}_j|0\rangle=\delta_{ij}
$$
\begin{eqnarray}
 a_i\bar{a}_{i_1}\cdots \bar{a}_{i_n}|0\rangle & = &\sum_{k=1}^{n}
\delta_{ii_k} \sum_{\pi \in S_{n,k}}
 \Phi^{i}_{i_1 \cdots i_n;\pi \cdot (i_1 \cdots \not \!{i_k} \cdots i_n)}\,
 \pi \cdot (\bar{a}_{i_1} \cdots \not \!{\bar{a}}_{i_k} \cdots
 \bar{a}_{i_n})|0\rangle \equiv \nonumber \\
&\equiv & \sum \tilde{\Phi}_{i_1 \cdots i_n;(n-1)}^{i}
\end{eqnarray}
$$
 \Phi^{i}_{j,0}=1
  $$
 where $\{\Phi \}$ are complex coefficients,$S_{n,k}$ denotes the group of
  permutations acting on $(i_1 \cdots \not \!{i_k} \cdots i_n)$ and
 the slash denotes the omission of the corresponding operator $\bar{a}_{i_k}$.
  The state
 in Eq.(3.15) is zero if the index ${\it i}$ is different from any of the
indices $(i_1 \cdots i_n)$.
 The matrix element of any monomial in $a_i$,$\bar{a}_i$ is equal to zero if
the
 indices in the monomial do not appear in pairs.For example,
 \begin{equation}
 \begin{array}{c}
 A_{i_1\cdots i_n;j_1 \cdots j_n}=\sum_{(n-1)}\sum_{(n-2)} \cdots \sum_1
 \tilde{\Phi}^{i_1}_{(i_1 \cdots
i_n),(n-1)}\tilde{\Phi}^{i_2}_{(n-1),(n-2)}\cdots
 \tilde{\Phi}^{i_n}_{(1),0}\equiv (\tilde{\Phi}^{(n)})!\\[4mm]
 \tilde{\Phi}^{i}_{j,0}=\delta_{ij}
 \end{array}
 \end{equation}
 is zero if $(j_1 \cdots j_n)$ is not a permutation of $(i_1 \cdots i_n)$.
 This is a necessary and sufficient condition for the existence of the
 number operators $N_i$,$\forall i \in I$,in a given Fock-like
representation.\\
 From Eq.(3.15),it is easy to see the action of any monomial $P(a_i,\bar{a}_j)$
on
 any state in a Fock-like space, particularly the action of the
operators$^{15}$
 $\Gamma_{ij}=a_i\bar{a}_j$.
 Hence,using Eq.(3.15) we find the relation expressing $a_i\bar{a}_j$ in terms
of normal
 ordering of operators a and $\bar{a}$:
 $$
\Gamma_{ij} \equiv a_i\bar{a_j}=\delta_{ij} + C^{ij}\bar{a}_ja_i
+ C^{ij}_{jk,ki}\bar{a}_j\bar{a}_ka_ka_i+
 C^{ij}_{jk,ik}\bar{a}_j\bar{a}_ka_ia_k
 $$
 $$
 +C^{ij}_{kj,ik}\bar{a}_k\bar{a}_ja_ia_k+
 C^{ij}_{kj,ki}\bar{a}_k\bar{a}_ja_ka_i
 + \cdots
$$
\begin{equation}
\equiv  \delta_{ij} + \somk \sum_{I^{n-1}}C^{ij}_{(\cdots j \cdots),(\cdots i
\cdots)}
(\bar{a}\cdots \bar{a}_j\cdots\bar{a})_k(a\cdots a_i \cdots a)_k,
\end{equation}
where any allowed monomial appearing on the RHS has the $\bar{a}_ja_i$
structure while
the rest of operators appear in  pairs of the same indices.Equations(3.15) and
Eq.(3.17)
are completely equivalent.They represent the most general relations between the
operators
$a_i$, $\bar{a}_j$ compatible with the number operators $N_i$,Eq.(3.6).
Let us remark that using Eqs.(3.6) and (3.15) we can express the number
operators
as

\begin{equation}
N_i=\bar{a}_ia_i + \sum_{k=2}^{\infty} \sum_{I^{n-1}} D^{i}_{(\cdots i \cdots),
(\cdots i \cdots)}
(\bar{a}\cdots \bar{a}_i \cdots \bar{a})_k \;
(a\cdots a_i \cdots a)_k
\end{equation}

and the transition number operators as
\begin{equation}
\begin{array}{c}
[N_{ij},\bar{a}_k]=\delta_{jk}\bar{a_i},\quad

[\tilde{N}_{ij},a_k ]= -\delta_{ik}a_j,\quad

{\tilde{N}}_{ij}=N_{ji}^{\dagger},\\[4mm]

N_{ij}=\bar{a_i}a_j + \sum_{k=2}^{\infty} \sum_{I^{n-1}} D^{ij}_{(\cdots i
\cdots),(\cdots j \cdots)}
(\bar{a}\cdots \bar{a}_i \cdots \bar{a})_k\, (a\cdots a_j \cdots a)_k,\\[4mm]

\tilde{N}_{ij}=\bar{a}_ia_j +\sum_{k=2}^{\infty}  \sum_{I^{n-1}}
\tilde{D}^{ij}_{(\cdots i \cdots),
(\cdots j \cdots)}
(\bar{a}\cdots \bar{a}_i \cdots \bar{a})_k \; (a\cdots a_j \cdots a)_k.\\
\end{array}
\end{equation}
{}From the above equation for $N_{ij}$ we can obtain Eqs.(3.15) and
(3.17).Hence,
 the sets of quantities $\{\Phi \}$,$\{ C \}$ and $\{ D \}$ in
Eqs.(3.15),(3.17)
 ,(3.18)
 are completely equivalent,as the set of quantities
 $\{\varphi_i \}$,$\{ c_i \}$ and $\{ d_i \}$ are equivalent for the
single-mode
 oscillator. Note that the relation $\bar{a}a=\varphi(N)$ has no simple
generalization
  for the multimode case.\\
In this paper we restrict ourselves to those relations between  the operators
$a_i$, $\bar{a}_i$ compatible with the number operators $N_i$,Eqs.(3.6).
However there remains an interesting question,namely  under which conditions
the set of
operators $a_i$, $\bar{a}_i$ with the most general contractions ,Eq.(3.15),
without number operators,can be obtained from the operators $c_i$, $\bar{c}_i$
with
number operators (e.g. by mapping).

In order that the constructed Fock-like space should become physically
acceptable,we
have to define the norm of vectors and demand that any linear combination of
vectors should have  non-negative norms.Therefore ,we demand a well-defined
scalar product on Fock-like space.This will generally lead to additional
restrictions
on the coefficients  $\{\Phi \}$,$\{ C \}$ and $\{ D \}$ .In principle,we have
two
possibilities.

We can demand that $\bar{a}_i$ should become Hermitian conjugate of the $a_i$
operator,with
respect to the scalar product $<*,*>$,i.e.$\bar{a}_i=a^{\dagger}_i$.In this
case the block matrix ${\bf A}$ becomes
Hermitian,with non-negative diagonal matrix elements,i.e.
\begin{equation}
A_{(i),(j)}=A_{(j),(i)}^{*} \qquad A_{(i),(i)}\geq 0.
\end{equation}
We define the norm  $||v||$ of any vector $|v>$ as
\begin{equation}
||v||^2=<v|v> \;  >0 \qquad  \;|v>\neq 0.
\end{equation}
Hence,the non-negative norms of all vectors leads us to the semi-positivity
   of
the matrix ${\bf A}$ in Eq.(3.20).This means that the matrix ${\bf A}$ has no
negative
eigenvalues.Particularly,any zero eigenvalue implies a null-vector (the
state of zero norm) which implies a relation between monomials in the operators
$\bar{a}_i$.\\
A general null-vector ${\bf E}$ is of the form
\begin{equation}
\sum_{\sigma \in S_n/S_t} A_{i_1 \cdots i_n,\sigma \cdot (i_1 \cdots i_n)}\,
E_{i_1 \cdots i_n,\sigma \cdot (i_1 \cdots i_n)}=0.
\end{equation}

Then ,the consistency conditions for null-vectors are
\begin{equation}
\sum_{\sigma \in S_n/S_t}E_{i_1 \cdots i_n,\sigma \cdot (i_1 \cdots i_n)}\;
\tilde{\Phi}^{j}_{\sigma \cdot (i_1 \cdots i_n),j_1 \cdots j_{n-1}}=0,\qquad
\forall
j \in I.
\end{equation}
This  means that any contraction of the operator $a_j$ with a null-vector is a
null-vector.
Hence,starting with   general contractions ,Eq.(3.15),with
$\bar{a}_i=a_i^{\dagger}$
leading to the Hermitian matrix ${\bf A}$ and positive diagonal elements
(Eq.(3.20)),
the positivity condition and consistency of null-vectors imply a strong
restriction
on the possible coefficients $\{\Phi \}$ (or $\{ C \}$ and $\{ D \}$).Eq.(3.15)
(or
(3.17) or (3.18)) uniquely determines the matrix ${\bf A}$,Eq.(3.20).However,
 the reverse is not generally true,
namely,even if the matrix ${\bf A}$ is Hermitian and positive definite  this
does not
mean that the corresponding contractions (3.15) exist.An interesting question
is whether
the above conditions for  the matrix ${\bf A}$,Eqs.(3.20) and (3.23),are also
sufficient
 for the existence of relations between $a_i$ and $a_i^{\dagger}$ of the type
 given in Eq.(3.15) or Eq.(3.17).\\
 The second interesting possibility is  that $\bar{a}_i \neq a_i^{\dagger}$,$i
\in I$.
 In this case there is no condition on the matrix ${\bf A}$ to be Hermitian and
positive
 definite.Moreover,the matrix ${\bf A}$ can have complex eigenvalues.Then,we
perform
 a polar decomposition of the matrix ${\bf A}$:
 $$
 A=U\cdot H
$$
where
$$
U\cdot U^{\dagger}=U^{\dagger}\cdot U=1,\qquad A\cdot A^{\dagger}=H^2 \geq 0.
$$
${\bf U}$ is a unitary matrix and
 ${\bf H}$ is a unique Hermitian,positive definite (non-negative) matrix.
Now the problem is to find out in which cases the
matrix ${\bf A}$ can lead to a relation between $a_i$ and $a_{i}^{\dagger}$
 of the type
 given by  Eq.(3.15) or Eq.(3.17).
 In the following we restrict ourselves to the case $\bar{a}_i=a_i^{\dagger}$ .

Finally,let us briefly describe the connection of the spectrum of the matrix
${\bf A}$ (Eq.(3.11)) with   generalized statistics,i.e. with the counting of
allowed
multiparticle states.As we have already stated,the appearance of null-vectors
implies
 corresponding relations between monomials in ${a_i}^{\dagger}$.If we fix the
indices
 $i_1\leq i_2 \cdots \leq i_n \in I$ with multiplicities $m_1$,$m_2 \cdots m_s$
 such that $\sum_{k=1}^{s} m_k=n$,then  the
 number of   linear independent states $\pi \cdot  (a_{i_1}^{\dagger} \cdots
a_{i_n}^{\dagger}) |0\rangle $ , $\pi \in S_n/S_t$, is equal to
\begin{equation}
 rank [A_{(i_1 \cdots i_n) fixed}]\leq \frac{n!}{m_1!m_2! \cdots m_s!}.
\end{equation}
 If we restrict the indices $i_1$,$i_2 \cdots i_n \in I_d \subset I$,where
${\it d}$ denotes
 the number of   single-mode oscillators,then the total number of n-particle
 excited states of d-oscillators is given by
 \begin{equation}
 W(n,d,A)=rank [A_{(n,d) }]\leq d^n.
 \end{equation}
 For example,for d-Bose oscillators
 \begin{equation}
 rank [A_{(i_1 \cdots i_n) fixed}]=1 ,\qquad
W_{B}(n,d)=\frac{(d+n-1)!}{n!(d-1)!}.
 \end{equation}
 For d-Fermi oscillators
 \begin{equation}
 rank [A_{(i_1 \cdots i_n) fixed}]=\theta(2-m_1)\cdots \theta(2-m_s),
 \qquad W_F(n,d)=\frac{d!}{n!(d-n)!}
 \end{equation}
 where $\theta (x)$ is the step-function, i.e. $\theta(x)=1$,$x>0$ and
$\theta(x)=0$
 $x \leq 0$.\\
 For ${\it n}$ quonic oscillators$^{13,14}$ which satisfy the relation $a_i
a_{j}^{\dagger}-
 q_{ij}a_{j}^{\dagger}
 a_i=\delta_{ij}$,$|q_{ij}| <1$,$q_{ij} \in {\bf C}$
\begin{equation}
  rank [A_{(i_1 \cdots i_n) fixed}]=\frac{n!}{m_1!m_2! \cdots m_s!}, \qquad
  W_Q(n,d)=d^n.
 \end{equation}

 %%%%%%%%%%%%%%%%%%%%%%%%%%%%%%%%%%%%%%%%%%%%%%%

\section{Algebras from Generalized Jordan-Wigner Transformations}
\setcounter{equation}{0}

We assume that the number operators $N_i$,Eq.(3.6), exist and that
\begin{equation}
a_ia^{\dagger}_j=q_{ij}\,a^{\dagger}_ja_i \qquad i\neq j
\qquad q^{*}_{ij}=q_{ji}
\end{equation}
Furthermore,we assume that the algebra (4.1) can be obtained
  by mapping from  the Bose algebra
\begin{equation}
 \begin{array}{c}
 [b_{i},b_{j}^{+}]=\delta_{ij},\quad \forall i,j \in I\\[4mm]
[b_{i},b_{j}]=0.
\end{array}
\end{equation}
Then,it follows that
\begin{equation}
a_{i}=b_{i} e^{\sum_{j}c_{ij}N_{j}}\sqrt{\frac {\varphi_{i}(N_{i})}{N_{i}}},
\end{equation}
where $c_{ij}$ are complex numbers and $\varphi_{i}(N_{i})$ are arbitrary
(complex) functions with the properties $\varphi_{i}(0)=0$,
$\lim_{\epsilon \rightarrow 0}\frac
{\varphi_{i}(\epsilon)}{\epsilon}=1$, $|\varphi_{i}(1)|=1, \forall i \in I.$
It is important to note that the number operators are preserved, i.e.
\begin{equation}
N_{i}^{(a)}=N_{i}^{(b)} \equiv N_{i}, \quad \forall i \in I.
\end{equation}
Then it is easy to find the corresponding deformed-oscillator algebra:
\begin{equation}
\begin{array}{c}
a_{i}a_{j}=e^{c_{ji}-c_{ij}} a_{j}a_{i}  \quad i \neq j\\[4mm]
a_{i}a^{+}_{j}=e^{c_{ij}+c_{ji}^{*}}a^{+}_{j}a_{i},  \quad i \neq j \\[4mm]
a_{i}a^{+}_{i}=|\varphi_{i}(N_{i}+1)|  e^{\sum_{j}(c_{ij} +c_{ij}^{*})N_{j}}
e^{(c_{ii}+c_{ii}^{*})}
\\[4mm]
a^{+}_{i}a_{i}=|\varphi_{i}(N_{i})|  e^{\sum_{j}(c_{ij}+c_{ij}^{*})N_{j}}.
\end{array}
\end{equation}
Generally, there are other mappings of the Bose algebra (Eqs. (4.2) ),
 but, in general, they do not have the number operators $ N_{i}^{(a)},$ and
Eq. (4.4) does not hold for mappings other than those in Eq.  (4.3).

We point out that the complete deformed-oscillator algebra is associative owing
to the mapping of the Bose algebra. The Fock space for the deformed-oscillator
algebra is spanned by powers of the creation operators $ a_{i}^{+}$, $ i
\in I$,
acting on the vacuum $|0>^{(a)}=|0>^{(b)} \equiv |0>.$ The states in the Fock
space are specified by the eigenvalues of the number operators $N_{i}$, namely
$|n_{1},n_{2},...n_{i}....>^{(a)}=|n_{1},n_{2},...n_{i}....>^{(b)}$.
(If there exists a number
$n_{i}^{(0)} \in {\bf N}$, such that $\varphi_{i}(n_{i}^{(0)})=0$, then
$N_{i}=0,1,...
(n_{i}^{(0)}-1)$.) \\
The states with unit norm are
$$
|n_{1}.....n_{n}>=\frac{(a_{1}^{+})^{n_{1}}.......(a_{n}^{+})^{n_{n}}}{\sqrt
{[\tilde{
\varphi}_{1}(n_{1})]!.....[\tilde{\varphi}_{n}(n_{n})]! }}
e^{-\frac{1}{2}\sum_{j}\theta_{ji}(c_{ij}+c_{ij}^{*})n_{i}n_{j}}|0,0,....0>\\[4mm]
$$
\begin{equation}
\begin{array}{c}
[\tilde{\varphi}(n)]! =\tilde{\varphi}(n)
\tilde{\varphi}(n-1)......\tilde{\varphi}(1)\\[4mm]
\tilde{\varphi}_{i}(n_{i})=|\varphi_{i}(n_{i})|
e^{(c_{ii}+c_{ii}^{*})n_{i}},
\end{array}
\end{equation}
where $\theta_{ij}$ is the step function.
(For anyons in $(2+1)$-dimensional space, $\theta $ is the angle function.)

Furthermore, the matrix elements of the operators $a_{i}, a_{i}^{+}$, $i \in
I$,
are
\begin{equation}
\begin{array}{c}
<....(n_{i}-1)....|a_{i}|....n_{i}...>=<...n_{i}....|a_{i}^{+}|....(n_{i}-1)....>^{*}\\[4mm]
=\sqrt{\varphi_{i}(n_{i})}e^{\frac{1}{2}\sum_{j}(c_{ij}+c_{ij}^{*})n_{j}}
\prod_{j\neq i}
\delta_{n_{j},n_{j}^{'}}.
\end{array}
\end{equation}
 For any $k=0,1,2,..$,  we also find that
$$
(a_{j}^{+})^{k}
(a_{j})^{k}=\frac{[\tilde{\varphi}_{j}(N_{j})]!}{[\tilde{\varphi}_{j}(N_{j}-k)]!}
e^{k\sum_{l\neq j}(c_{jl}+c_{jl}^{*})N_{l}}\\[4mm]
$$
\begin{equation}
(a_{j})^{k}
(a_{j}^{+})^{k}=\frac{[\tilde{\varphi}_{j}(N_{j}+k)]!}{[\tilde{\varphi}_{j}(N_{j})]!}
 e^{k\sum_{l\neq j}(c_{jl}+c_{jl}^{*})N_{l}} .
\end{equation}
The norms of arbitrary linear combinations of the states in Eq. (4.6)
in the Fock space, corresponding to the deformed-oscillator algebra,
 are positive by definition
owing to the mapping of the Bose algebra (Eqs. (4.2) and (4.3)). \\Namely,
$|n_{1},n_{2},...n_{i}....>^{(a)}=
=|n_{1},n_{2},...n_{i}....>^{(b)} \equiv |n_{1},n_{2},...n_{i}....>$.

This class of deformed multimode  oscillator algebras comprises
multimode\\ Biedenharn-Macfarlaine , Aric-Coon , two-($p,q$)
parameter , Fermi,genons,
generalized Green's , as well as anyonic  and Pusz-Woronowicz
 oscillators  covariant under the $SU_{q}(n)$
$(SU_{q}(n|m))$ algebra (superalgebra)$^{4,7,8,9}$ .Particularly,the
operator algebra for the Haldane exclusion statistics$^{18}$ is a special
case of our mapping,Eqs.(4.3) and (4.5),with the substitutions

$$
c_{ij}=(-)c_{ji}=(-) \frac{\imath \pi}{m+1} \quad i<j, m \in {\bf N};
$$

\begin{equation}
K(N_i,g)=\frac{\varphi (N_i+1,g)}{\varphi (N_i,g)} \qquad g=\frac{1}{m}.
\end{equation}

Nonisomorphic (nonequivalent) algebras are classified by different matrix
elements given by (4.7), i.e. with the functions $ g_{i}(n_{1},n_{2},...n_{n})=
|\varphi_{i}(n_{i})|e^{\sum_{j}(c_{ij}+c_{ij}^{*})n_{j}}.$ It is important to
mention that there are mappings of the Bose algebra which do not preserve the
relation
$N_{i}^{(a)}=N_{i}^{(b)}$,  given by (4). Moreover, there are mappings of Bose
algebra
for which the number operators $N_{i}^{(a)} $ do not   exist. Such an example
is the exchange algebra presented in Ref. (19)

%%%%%%%%%%%%%%%%%%%%%%%%%%%%%%%%%%%%%%%%%%%%%%%%%%%%%%%%%%%%%%%%%

\section{Parabosons and parafermions}
\setcounter{equation}{0}

  The para-algebra is defined by the trilinear commutation relation$^{11}$
\begin{equation}
\begin{array}{c}
[a_ia^{\dagger}_j+q\,a^{\dagger}_j a_i,a_k]_{-}=-\frac{2}{p}q\delta_{jk}a_i
\\[4mm]

a_ia^{\dagger}_{i_1}a^{\dagger}_{i_2}\cdots a^{\dagger}_{i_n}|0\rangle=
\sum_{k=1}^{n} \delta_{ii_k}(-q)^{k-1}a^{\dagger}_{i_1}\cdots
\hat{a}^{\dagger}_{i_k}\cdots
 a^{\dagger}_{i_n}|0\rangle -\\[4mm]

 -\frac{2}{p}\sum_{k=2}^{n}\delta_{ii_k}\sum_{l=1}^{k-1}
 (-q)^la^{\dagger}_{i_1}\cdots
a^{\dagger}_{i_{k-1}}a^{\dagger}_{i_l}a^{\dagger}_{i_{k+1}}\cdots
a^{\dagger}_{i_n}|0\rangle.
\end{array}
\end{equation}
Here,$p>0$ is the order of parastatistics, $q=-1$ for the para-Fermi algebra
and $q=+1$
for the para-Bose algebra.\\
Comparing with Eq.(3.15),one  identifies the set $\{\Phi\}$ as
\begin{equation}
\begin{array}{c}
\Phi^{i}_{i_1\cdots i_n,id\cdot(i_1\cdots \not \!i_k \cdots i_n)}=(-q)^{k-1}
 \qquad 1\leq
k\leq n ,\\[4mm]

\Phi^{i}_{i_1\cdots i_n,\pi_{l,k-1}\cdot
(i_1\cdots \not \! i_k \cdots i_n)}=-(\frac{2}{p})(-q)^l \qquad 1\leq l\leq
k-1,
\end{array}
\end{equation}
where
$$
 \pi_{l,k-1}=  \left ( \begin{array}{ccc}
l &k-2&k-1\\
l+1&k-1&l\\
\end{array} \right )
$$

denotes the cyclic permutation and ${\it id}$ denotes the identity
permutation.\\
The matrix ${\bf A}$ ,Eq.(3.11),for the three different para-oscillators is

$$
A=  \left ( \begin{array}{cccccc}
1 & x & x & x^2 & x^2 &z\\
x & 1 & x^2 & x & z & x^2\\
x & x^2 & 1 & z & x & x^2\\
x^2 & x & z & 1 & x^2 & x\\
x^2 & z & x & x^2 & 1 & x\\
z & x^2 & x^2 & x & x & 1
\end{array} \right )
$$

\begin{equation}
\begin{array}{c}
x=q\,(\frac{2}{p}-1) \equiv qy ,\\[4mm]
 z=q\,\frac{2}{p}-q^3\,(\frac{2}{p}-1)^2.
\end{array}
\end{equation}
The matrix is written in the basis
$a^{\dagger}_1a^{\dagger}_2a^{\dagger}_3|0\rangle$,
$a^{\dagger}_2a^{\dagger}_1a^{\dagger}_3|0\rangle$,
$a^{\dagger}_1a^{\dagger}_3a^{\dagger}_2|0\rangle$,
$a^{\dagger}_2a^{\dagger}_3a^{\dagger}_1|0\rangle$,
$a^{\dagger}_3a^{\dagger}_1a^{\dagger}_2|0\rangle$,
 and $a^{\dagger}_3a^{\dagger}_2a^{\dagger}_1|0\rangle$.

By inspection of the eigensystem of the matrix ${\bf A}$,Eq.(5.3),one finds
 that the
rank ${\bf A} =4$
for $q=\pm 1$ and $\forall p \in {\bf N}$, which means that only four of the
states
  are linearly
independent.This is in accordance with the trilinear relations which hold for
 parastatistics $^{11}$,namely
$ [a_k,[a_m,a_n]_{\pm}]=0$.
 We choose
 the  set $a^{\dagger}_1a^{\dagger}_2a^{\dagger}_3|0\rangle$,
 $a^{\dagger}_2a^{\dagger}_1a^{\dagger}_3|0\rangle$,
$a^{\dagger}_1a^{\dagger}_3a^{\dagger}_2|0\rangle$ and
$a^{\dagger}_3a^{\dagger}_2a^{\dagger}_1|0\rangle$ as linearly
independent vectors. Instead of the $6\times 6$
matrix ${\bf A}$
 ,in the following we   use the $4\times 4$ matrix,corresponding to this set
  of vectors.\\ In the limit $p\rightarrow \infty$,the rank of ${\bf
A}$,Eq.(5.3),
  reduces to 1 for $q=\pm 1$, i.e. the matrix ${\bf A}$ reduces to Fermi (Bose)
  matrix for $q=+1$ (q=-1).Particularly,for the two-level system
$a^{\dagger}_i$,$i=1,2$,
  Ref.(12), the matrix ${\bf A}$,corresponding to the
 $a^{\dagger}_1a^{\dagger}_2|0\rangle$ and
  $a^{\dagger}_2a^{\dagger}_1|0\rangle$ states is

  $$
  A= \left ( \begin{array}{cc}
  1 & x\\
  x & 1
  \end{array} \right )
  $$

 For $q=\pm 1$, rank ${\bf A}=2 $,$ \forall p >1$,
   and again,in the limit $p\rightarrow \infty$ the matrix ${\bf A}$ is of the
Fermi (Bose) type for $q=+1$ (q=-1).\\
For the states $(a^{\dagger}_1)^2a^{\dagger}_2|0\rangle$ and
$ a^{\dagger}_1a^{\dagger}_2 a^{\dagger}_1 |0\rangle$,the matrix ${\bf A}$ is

$$
 A= \left ( \begin{array}{cc}
 1+x & x+x^2\\
 x+x^2 & 1+z
 \end{array} \right )
 $$
\baselineskip=24pt

For $q=+1$,rank ${\bf A}=2$,$\forall p > 1$. For $q=-1$,rank ${\bf A}=2$,
 p $\geq $ 3 and rank ${\bf A}=1$ if $p=1,2$.For para-Fermi case,when $p=2$
there exist
 additional relation between operators,
 what was overlooked in the paper Ref.(12).Similar statements hold for the
 matrix ${\bf A}$ corresponding to the states
  $(a^{\dagger}_1)^r(a^{\dagger}_2)^s|0\rangle$,$r+s=n$.As a consequence,the
para-Fermi
   (para-Bose) system behaves as Bose (Fermi) system when $p\rightarrow \infty$
$^{10}$.\\
Next,we discuss the number operators for para-Bose and para-Fermi
oscillators.\\
Number operators $N_i$,$i \in I$,can be written as in Eq.(3.18).We have,up to
the
second order
in the creation and annihilation operators ($k=2$ in Eq.(3.18)),
\begin{equation}
N_i=a^{\dagger}_ia_i +\frac{p^2}{4(p-1)} \sum_{l}[Y_{il}]^{\dagger}[Y_{il}],
\end{equation}
where
$ Y_{il}=a_ia_l-x\,a_la_i$.\\In the second order,we identify the set $\{ D \}$
as
$$
D^i_{li,il}=\frac{p^2}{4(p-1)};\qquad  D^i_{il,il}=D^i_{li,li}=q\frac{p(p-2)}
{4(p-1)};
\qquad  D^i_{il,li}=\frac{(p-2)^2}{4(p-1)}.
$$

The computation of the terms of the higher order ($k\geq 3$) becomes more
involved.For example,the generic form of the third order term is
$$\sum_{j,l}\sum_{\pi} X_{\pi\cdot(ijl)}\,\pi \cdot(a_la_ja_i).$$For the
parasystem,
the summation is taken over those set of indices which are linearly
independent,
owing to the constraint implied by trilinear relation\\ $
[a_k,[a_m,a_n]_{\pm}]=0$.
\\We illustrate the procedure by calculating the third order  term for the
$N_1$,
i.e. for $\{i,j,l\}=\{1,2,3\}$.The basic set of the X-operators is
$\hat{X}\equiv \{ X_{123},X_{213},X_{132},X_{321}\}$,which is solutions
of the equation

\begin{equation}
A\cdot \hat{X}=V
\end{equation}
where
$$
V=  \left ( \begin{array}{c}
0\\
-q(\frac{2}{p})v_1-(\frac{2}{p})v_2+(\frac{2}{p})(\frac{2}{p}-1)v_3+(\frac{2}{p})v_6\\
0\\
-q(\frac{2}{p})^2 v_1+(\frac{2}{p})(\frac{2}{p}-1)v_3+(\frac{2}{p})v_6
\end{array} \right )
$$

$$
q=\left \{ \begin{array}{cc}
  +1 & \mbox{para-Bose}\\
  -1 &\mbox{para-Fermi}
\end{array}
\right.
$$

\begin{equation}
v_1= a^{\dagger}_1a^{\dagger}_2a^{\dagger}_3 \quad
v_2=a^{\dagger}_2a^{\dagger}_1a^{\dagger}_3 \quad
v_3=a^{\dagger}_1a^{\dagger}_3a^{\dagger}_2\quad
v_6=a^{\dagger}_3a^{\dagger}_2a^{\dagger}_1
\end{equation}
For the parabosons (q=+1;$y=(\frac{2}{p}-1))$ one finds
\begin{eqnarray}
X_{123}&=&\frac{1}{2y-y^2-y^3}[-(1+3y+y^2+y^3)v_1+(y^2-y)v_2+(y+2y^2)v_3+(1+2y)v_6] \nonumber\\[6mm]
X_{213}&=&\frac{1}{-2+y+y^2}[(1-y)v_1+2v_2-yv_3-v_6]\nonumber \\[6mm]
X_{132}&=&\frac{1}{-2+y+y^2}[-(1+2y)v_1-yv_2+(y+y^2)v_3+(1+y)v_6]\nonumber
\\[6mm]
X_{321}&=&\frac{1}{y}X_{132}
\end{eqnarray}
The corresponding operators $\hat{X}$ for the parafermions are obtained from
the above equations by substitutions $ X_{123}\Rightarrow (-)X_{123}(-y,-v_1)$,
$X_{213}\Rightarrow X_{213}(-y,-v_1)$,$X_{132}\Rightarrow X_{132} (-y,-v_1)$,
$X_{321}\Rightarrow X_{321}(-y,-v_1)$.\\
The operators $\hat{X}$ with some of the indices equal (e.g.
$X_{122}$,$X_{212}$,$X_{111}$ etc.)
follow from the Eq.(5.7) after the suitable identification of the indices.

The transition number operators $N_{ij}$,$i,j \in I$, are given,up to the
second order,by
\begin{equation}
N_{ij}=a^{\dagger}_ia_j +\frac{p^2}{4(p-1)} \sum_{l}[Y_{il}]^{\dagger}[Y_{jl}]
\end{equation}
To this order,we identify the set $\{ D \}$
as
$$
D^{ij}_{li,jl}=\frac{p^2}{4(p-1)}\qquad  D^{ij}_{il,jl}=D^i_{li,lj}=
q\frac{p(p-2)}{4(p-1)}
\qquad  D^i_{il,lj}=\frac{(p-2)^2}{4(p-1)}
$$
Similarly as in the case of the number operator $N_1$,we compute
the basic set of the third order operators $\tilde{X}=\{ \tilde{X}_{123},
\tilde{X}_{213}, \tilde{X}_{132},\tilde{X}_{321} \}$ for the transition
number operator $N_{12}$.These operators are the solutions of the equation
$$
A\cdot \tilde{X}=\tilde{V}
$$
where
$$
\tilde{V}=\left ( \begin{array}{c}
-(q+x)w_1+x(q+x)w_2\\
0\\
x(q+x-1)\,w_1+(1-q-x-qx^2)\,w_2+qx\,w_3\\
(x-1)(q+x)\,w_1
\end{array} \right )
$$

\begin{equation}
\begin{array}{c}
w_1=a^{\dagger}_1a^{\dagger}_1a^{\dagger}_3 \\
w_2=a^{\dagger}_1a^{\dagger}_3a^{\dagger}_1 \\
w_3=a^{\dagger}_3a^{\dagger}_1a^{\dagger}_1
\end{array}
\end{equation}

$$
w_3=\left \{\begin{array}{cc}
w_1& \mbox{para-Bose}\\
2w_2-w_1   &\mbox{para-Fermi}.
\end{array}
\right.
$$

For the parabosons ($x=qy=(\frac{2}{p}-1)$) the result is:
\begin{eqnarray}
\tilde{X}_{123}&=&\frac{1}{x^2+x-2}[-3x\,w_1+(1+x+x^2)\,w_2] \nonumber \\[6mm]
\tilde{X}_{213}&=&\frac{1}{x^2+x-2}[x^2 \,w_1-x\,w_2]\nonumber \\[6mm]
\tilde{X}_{132}&=&\frac{1}{x^2+x-2}[-2x\,w_1+(x+x^2)\,w_2]\nonumber \\[6mm]
\tilde{X}_{321}&=&(-)\frac{1+x}{x}\tilde{X}_{132}.
\end{eqnarray}

For the parafermions ($x=qy=-(\frac{2}{p}-1)$) the result is:
\begin{eqnarray}
\tilde{X}_{123}&=&\frac{1}{x^3-x^2-2x}[(-5x^2+2x-2)\,w_1
+(x^3+x^2+3x)\,w_2]\nonumber\\[6mm]
\tilde{X}_{213}&=&\frac{1}{x^2-x-2}[(x^2+2)\,w_1-3x\,w_2]\nonumber\\[6mm]
\tilde{X}_{132}&=&\frac{1}{x^2-x-2}[(2x-2)\,w_1+(-x^2+x-4)\,w_2]\nonumber\\[6mm]
\tilde{X}_{321}&=&\frac{1}{x^3-x^2-2x}[(-2x^2+2x-2)\,w_1+(x^3-2x^2+3x)\,w_2]
\end{eqnarray}

Again,the operators $\tilde{X}$ with some of the indices equal (e.g.
$\tilde {X}_{122}$,$\tilde {X}_{212}$,$\tilde {X}_{111}$ etc.)
follow from the above result after the suitable identification of the indices.

%%%%%%%%%%%%%%%%%%%%%%%%%%%%%%%%%%%%%%%%%%%%%%%%%%%%%%%%%%%%%%%%%%%
%%%%%%%%%%%%%%%%%%%%%%%%%%%%%%%%%%%%%%%%%%%%%%%%%%%%%%%%%%%%%%%%%%%%

\section{ Paraquantization of Govorkov}
\setcounter{equation}{0}
\baselineskip=24pt

Govorkov$^{13}$ has defined a new para-algebra obeying

\begin{equation}
[a_{i} a^{\dagger}_{j}, a_{k}] = (\frac {\lambda}{p}) \delta_{jk} a_{i}
\end{equation}
$$
a_{i}a^{\dagger}_{i_1} a^{\dagger}_{i_2} \cdots a^{\dagger}_{i_n} |0 > =
\delta_{ii_1} a^{\dagger}_{i_2}\cdots a^{\dagger}_{i_n}| 0 > - (\frac
{\lambda}{p})
\sum^{n}_{k=2} \delta_{ii_k} a^{\dagger}_{i_2} a^{\dagger}_{i_3} \cdots
a^{\dagger}_{i_{k-1}} a^{\dagger}_{i_1} a^{\dagger}_{i_{k-1}} \cdots
 a^{\dagger}_{i_{k+1}}\cdots
a^{\dagger}_{i_n} |0>.
$$

Comparing with Eq. (3.15), one identifies the set ${\Phi}$ as

$$
\Phi^{i}_{i_{1}... i_{n}, id \cdot (i_{1}...\not \! i_{k} ... i_{n})} =
\delta_{k1}
$$
\begin{equation}
\Phi^{i}_{i_{1}...i_{n}, \pi_{l,k-1}\cdot (i_{1}...\not \! i_{k}...i_{n})}=
-(\frac {\lambda}{p}) \delta_{l1} \ \ \ 2 \leq k \leq n.
\end{equation}

With restriction to the three oscillators $(i_{k}=1, 2, 3)$, one immediately
obtains
the matrix {\bf A} as
\begin{displaymath}
A= \left (
\begin{array}{rrrrrr}
1 & -y & -y & y^{2} & y^{2} & -y\\
-y & 1 & y^{2} & -y & -y & y^{2}\\
-y & y^{2} & 1 & z & y & y^{2}\\
y^{2} & -y & -y & 1 & y^{2} & -y\\
y^{2} & -y & -y & y^{2} & 1 & -y\\
-y & y^{2} & y^{2} &-y & -y & 1\\
\end{array}
\right)
\end{displaymath}

where $y=\frac {\lambda}{p}$. The matrix $A$ is written in the same basis as
the
matrix $A$ for the para-Bose (para-Fermi) oscillators in Eq. (5.3).\\
The number operators $N_i$ can be written ,up to the second order
($k=2$ in Eq.(3.18)) as

\begin{equation}
N_i=a^{\dagger}_ia_i+\frac{1}{(1+\frac{\lambda}{p})(1-\frac{\lambda}{p})}
 \sum_{l}[Y_{il}]^{\dagger}[Y_{il}]
\end{equation}
where $Y_{ik}=a_{i} a_{k} + (\frac {\lambda}{p}) a_{k} a_{i}$.\\
To this order,we identify the set $\{ D \}$
as
$$
D^{i}_{li,il}=\frac{1}{(1+\frac{\lambda}{p})(1-\frac{\lambda}{p})};\\
 D^{i}_{il,il}=D^i_{li,li}=(\frac{\lambda}{p})\frac{1}
 {(1+\frac{\lambda}{p})(1-\frac{\lambda}{p})};\\
D^{i}_{il,li}=(\frac{\lambda}{p})^2\frac{1}
{(1+\frac{\lambda}{p})(1-\frac{\lambda}{p})}.
$$
For the number operator $N_1$ ($\{i,j,l\}=\{1,2,3\}$), the basic set of the
 third-order
operators $X_{\pi\cdot (ijl)}$, $\{ ijl \}=\{ 1,2,3 \}$, is
$\hat{X}\equiv \{ X_{123},X_{213},X_{132},X_{231}, X_{312},X_{321}\}$,which are
the solutions
of the equation

\begin{equation}
A \cdot \hat {X} = V \ \
\end{equation}

\begin{displaymath}
V = \left (
\begin{array}{c}
0\\
(\frac {\lambda}{p})^{2} v_{3} + (\frac {\lambda}{p}) v_{5}\\
0\\
(\frac {\lambda}{p})^{2} v_{1} +(\frac {\lambda}{p})v_{2}+ (\frac
{\lambda}{p})v_{3}+v_{4}\\
(\frac {\lambda}{p})^{2} v_{1} + (\frac {\lambda}{p}) v_{2}\\
(\frac {\lambda}{p})^{2} v_{3} + (\frac {\lambda}{p}) v_{1} + (\frac
{\lambda}{p})
v_{5} + v_{6}\\
\end{array}
\right)
\end{displaymath}
and with $v_{i}$ 's defined as $v_{1}=a^{\dagger}_{1} a^{\dagger}_{2}
a^{\dagger}_{3},
v_{2}=a^{\dagger}_{2} a^{\dagger}_{1} a^{\dagger}_{3}, v_{3}= a^{\dagger}_{1}
a^{\dagger}_{3} a^{\dagger}_{2}, v_{4}=a^{\dagger}_{2} a^{\dagger}_{3}
a^{\dagger}_{1},
v_{5}= a^{\dagger}_{3} a^{\dagger}_{1} a^{\dagger}_{2}$, and
$v_{6}=a^{\dagger}_{3}
a^{\dagger}_{2} a^{\dagger}_{1}$.\\
Explicitly,
\begin{eqnarray}
X_{123}&=&\frac{y[(y+4y^3)v_1+4y^2v_2+4y^2v_3+2yv_4+2yv_5+v_6]}{1-5y^2+4y^4}\nonumber \\[6mm]
X_{213}&=&\frac{y[4y^2v_1+2yv_2+2yv_3+v-4+v-5+2yv_6]}{1-5y^2+4y^4}\nonumber
\\[6mm]
X_{231}&=&\frac{2y^2v_1+yv_2+yv_3+(1-2y^2)v_4+2y^2v_5+yv_6}{1-5y^2+4y^4}
\end{eqnarray}

The remaining three X's are obtained from the above equation by the
substitutions $v_2\leftrightarrow v_5$,$v_1\leftrightarrow v_3$ and
$v_4\leftrightarrow v_6$,i.e.$ X_{123}\rightarrow X_{132}$,
$X_{213}\rightarrow X_{312}$ and $X_{231}\rightarrow X_{321}$.
Notice that,for $p>2$, there are six linearly independent X's since there
are no trilinear constraint as in the case of the parastatistics.

Transition number operators $N_{ij}$ are,up to the second order,
\begin{equation}
N_{ij}=a^{\dagger}_ia_j+\frac{1}{(1+\frac{\lambda}{p})(1-\frac{\lambda}{p})}
 \sum_{l}[Y_{il}]^{\dagger}[Y_{jl}] ,
\end{equation}
and the coefficients $\{ D \}$
are
$$
D^{ij}_{li,jl}=\frac{1}{(1+\frac{\lambda}{p})(1-\frac{\lambda}{p})};\quad
 D^{ij}_{il,jl}=D^i_{li,lj}=(\frac{\lambda}{p})\frac{1}
 {(1+\frac{\lambda}{p})(1-\frac{\lambda}{p})};\\
D^{ij}_{il,lj}=(\frac{\lambda}{p})^2\frac{1}
{(1+\frac{\lambda}{p})(1-\frac{\lambda}{p})}
$$
The basic set of the third-order operators $\tilde{X}$ for the transition
number
operator $N_{12}$ are

$$
\tilde{X}_{123}=X_{213}(1\equiv 2) \qquad \tilde{X}_{213}=X_{123}(1\equiv 2)
\qquad
\tilde{X}_{132}=X_{231}(1\equiv 2)
$$
$$
\tilde{X}_{231}=X_{132}(1\equiv 2) \qquad \tilde{X}_{312}=X_{321}(1\equiv 2)
\qquad
\tilde{X}_{321}=X_{312}(1\equiv 2)
$$
Here,the abbreviations,e.g.  $X_{213}(1\equiv 2)$, mean that one has to
identify
indices $1$ and $2$ in  $v_i$'s such that $v_1=v_2=
a^{\dagger}_1a^{\dagger}_1a^{\dagger}_3$,$
v_3=v_4=a^{\dagger}_1a^{\dagger}_3a^{\dagger}_1$ and $v_5=v_6=
a^{\dagger}_3a^{\dagger}_1a^{\dagger}_1$ and than read off X's from Eq.(6.5).\\

%%%%%%%%%%%%%%%%%%%%%%%%%%%%%%%%%%%%%

\section{Generalized Quon Algebra}
\setcounter{equation}{0}
\baselineskip=24pt

The general (associative) quon algebra Ref.(16) is defined by:
\begin{equation}
\begin{array}{c}
a_i a_j^{\dagger} -q_{ij} a_j^{\dagger} a_i = \delta_{ij} ,  \qquad \forall i,j
\in I ,\\[6mm]
q_{ij}^{\ast } = q_{ji},\\
a_l(a_{i_1}^{\dagger} \cdots a_{i_n}^{\dagger})|0>=\sum_{\alpha =1}^n
q_{li_1}\cdots
q_{li_{\alpha -1}}a_{i_1}^{\dagger}\cdots \not \!a_{i_{\alpha}}^{\dagger}
\cdots
a_{i_n}^{\dagger}\delta_{li_{\alpha}}|0>.

\end{array}
\end{equation}
No commutation relation between
$a_i$ and $a_j$ exists if $|q_{ij}| \neq 1$ for  $\forall i,j $ $\in I$ .\\
One identifies the set $\Phi$ as
\begin{equation}
\Phi^i_{i_1 \cdots i_n,id\cdot (i_1 \cdots \not \! i_k \cdots i_n)} =
\prod_{k=1}^{\alpha -1}
q_{lk},\qquad 1\leq \alpha \leq n.
\end{equation}

The matrix ${\bf A}$ is hermitian and block-diagonal.
A generic block
$A^{(i_{1}\cdots i_{n})}$ is characterizied
by mutually different ordered indices $i_{1},\dots i_{n}\in I$ $(i_{1}<
i_{2}< \cdots < i_{n})$ from which all other blocks in the
n-particle sector can be obtained using a suitable specification. The
$A^{(i_{1},\dots
i_{n})}$ matrix is an
$n!\cdot n!$ matrix, whose diagonal matrix elements are equal to $1$.
The arbitrary matrix element $(\pi ,\sigma)$, i.e.,$i_{\pi (1)},\cdots
i_{\pi (n)}; i_{\sigma (1)}\cdots i_{\sigma (n)}$, where $\pi $ and
$\sigma $ are permutations acting on positions $1,2\dots n$
($\pi $ denotes the row and $\sigma $ the column of the matrix $A^{(i_{1},\dots
i_{n})}$)
is given by

\begin{equation}
A^{(i_{1},\dots ,i_{n})}_{\pi ,\sigma } =
\prod_{\alpha, \beta} q_{i_{\alpha }i_{\beta }}.
\end{equation}

Here the product is over all pairs $\alpha ,\beta =1,\dots ,n$
satisfying $\pi^{-1}(\alpha) < \pi^{-1}(\beta)$ and
$\sigma^{-1}(\alpha) > \sigma^{-1}(\beta)$.

The general simple structure of the number operator $N_k$ ( $|q_{ij}| < 1$,
$ \forall i,j \in S$) is

\begin{equation}
N_{k} = a_{k}^{\dagger}a_{k} + \sum^{\infty}_{n=1}\sum_{i_{1},\dots ,i_{n}}
\sum_{\pi \in S_{n}}[Y_{k,\pi(i_{1},\dots ,i_{n})}]^{\dagger}
Y_{k,i_{1},\dots ,i_{n}}(A^{(k,i_{1},\dots ,i_{n})})^{-1}_{k,i_{1},\dots
,i_{n};
k,\pi(i_{1},\dots ,i_{n})} ,
\end{equation}
where
$$
Y_{ki_1}  = a_k a_{i_1} - q_{i_1k} \; a_{i_1} a_k
$$
$$
Y_{k,i_1\cdots i_n}=Y_{k,i_1\cdots i_{n-1}}a_{i_n}-q_{i_nk}q_{i_ni_1}\cdots
q_{i_ni_{n-1}}a_{i_n}Y_{k,i_1\cdots i_{n-1}}
$$
and ${\bf A}$ is matrix given in Eq.(7.3).\\
Always when the parameters $q_{ij}$ tend to $1$, i.e., $q_{ij}
\rightarrow 1$, for $\forall i,j$, quons tend to a particular anyonic-type
statistics, the matrices $(A^{(i_{1},\dots i_{n})})^{-1}$ become singular and
the coefficients
in the number operator, Eq.(7.4), diverge. Nevertheless, the number operator
$N_{k}$, when acting on states, is well defined. Moreover, in the exact limit
it reduces to $N_{k} = a_{k}^{\dagger}a_{k}$, and additional relations between
the annihilation (creation) operators $a_{i},a_{j}$ ($a_{i}^{\dagger},a_{j}
^{\dagger}$) emerge. In this case, the corresponding particles are not
distinguishable, i.e., they are identical in the quantum-mechanical sense.
Interchanging them, we generally obtain a unit phase factor $e^{i\alpha}$
(typical of anyons).\\
The transition number operator $N_{ij}$ has also simple structure:
\begin{equation}
N_{ij}=a^{\dagger}_ia_j+\sum^{\infty}_{n=1}\sum_{i_{1},\dots ,i_{n}}
\sum_{\pi \in S_{n}}[Z^i_{ji_{1},\dots ,i_{n}}]^{\dagger}
Y_{j,\pi i_{1},\dots ,i_{n}}(A^{(k,i_{1},\dots ,i_{n})})^{-1}
_{j\pi(i_{1},\dots ,i_{n});ji_{1},\dots ,i_{n}},
\end{equation}
where
$$
Z^k_{ki}=a_la_i-q_{ik}a_ia_l \equiv (Y_{ki})_{\{a_k\rightarrow a_l\}}
$$
$$
Z^l_{ki_{1},\dots ,i_{n}}\equiv (Y_{ki_{1},\dots ,i_{n}})_{\{a_k
\rightarrow a_l\}}
$$

%%%%%%%%%%%%%%%%%%%%%%%%%%%%%%%%%%%%%%%%%%%%%%%%%%%%%%%%%%%%%%%%%%%%%%%%%

\section{Conclusion}
In this paper we have generalized our previous analysis of  single-mode
oscillators$^6$ to
  multimode oscillator algebras with Fock-like representations.We have extended
the
  quantization defined by $a_ia^{\dagger}_j-R_{ij,kl}
  a^{\dagger}_ka_l=\delta_{ij}$,$i,j \in I$ to include all normally ordered
terms in
  creation and annihilation operators.We have restricted ourselves to the
algebras
   with  well-defined
  number operators and transition number operators.In this way we have unified
  all approaches to  multimode
  oscillator algebras $^{7-10,13-17}$.The connection between these algebras and
generalized
  statistics has been pointed out.
%%%%%%%%%%%%%%%%%%%%%%%%%%%%%%%%%%%%%%%%%%%%%%%%%%%%%%%%%%%%
\newpage
{\bf Aknowledgment}

We thank D.Svrtan for helpful discussions.
This work was supported by the Scientific Fund of the Republic of Croatia.

%%%%%%%%%%%%%%%%%%%%%%%%%%%%%%%%%%%%%%%%%%%%%%%%%%%%%%%%%%%%%%%%%%%%%%%%%


\newpage
\begin{thebibliography}{99}

\bibitem{1.} M.Jimbo,Lett.Math.Phys.{\bf 10}, 63 (1985);
V.G.Drinfeld,{\it Quantum
Groups} (Proc.Int.Congr.of Math.,Berkeley,CA,1986).

\bibitem{2.} S.L.Woronowicz, {\it Comm.Math.Phys.} {\bf 111}, 613 (1987);
 Yu.I.Manin,{\it Quantum Groups and Non-Commutative Geometry}
   (Center de Recherches,University of Montreal,1988);
  J.Wess and B.Zumino, {\it Nucl.Phys.Suppl.}{\bf 318},302 (1990);
 A.Connes, {\it Noncommutative Differential Geometry}
 (Cambridge University Press,Cambridge,1994).

\bibitem{3.} For a recent reviews,see : C.Gomez,M.Ruiz-Altaba and G.Sierra,
{\it Quantum Groups in Two-Dimensional Physics} (Cambridge University Press,
Cambridge,1993); Z.Chang, {\it Quantum Groups and Quantum Symmetry },
 preprint IC/94/89 (unpublished).

\bibitem{4.} M.Arik and D.D.Coon,
  {\it J.Math.Phys.}{\bf 17}, 524 (1976); L.C.Biedenharn,
   {\it J. Phys.A :Math.Gen.}{\bf 22}, L873 (1989); A.J.Macfarlane,
  {\it J. Phys.A :Math.Gen.}{\bf 22}, L983 (1989);
   P.P.Kulish and E.V.Damaskinsky,
   {\it J. Phys.A :Math.Gen.} {\bf 23}, 415 (1990);
 L.de Falco et.al.,{\it Mod.Phys.Lett.}{\bf A9}, 3331 (1994); V.I.Man'ko
et.al.,
{\it Int.J.Mod.Phys.}{\bf A8}, 3577 (1993); T. Brzezinski,I.L.Egusquiza and
A.J. Macfarlane,
{\it Phys.Lett.}{\bf 311B}, 202 (1993); K.Odaka,T.Kishi and S.Kamefuchi,
   {\it J. Phys.A :Math.Gen.}{\bf 24}, L591 (1991); C.Chou,
{\it Mod.Phys.Lett.}{\bf A7}, 2685 (1992); N. Debergh,
{\it Mod.Phys.Lett.}{\bf A8}, 765 (1993); E.Celeghini,T.D.Palev and M.Tarlini,
{\it Mod.Phys.Lett.}{\bf B5}, 187 (1991); R. Floreanini and L.Vinet,
 {\it J. Phys.A :Math.Gen.}{\bf 23}, L1019 (1990);
S.Chaturvedi,V.Srinivasan and R.Jagannathan,
 {\it Mod.Phys.Letters}{\bf A8}, 3727 (1993) .

\bibitem{5.} D.Bonatsos and C.Daskaloyannis,
 {\it J. Phys.A :Math.Gen.}{\bf 26}, 1589 (1993) ;
{\it Phys.Lett.}{\bf 307B}, 100 (1993).

\bibitem{6.} S.Meljanac,M.Milekovi\c $\,$and S.Pallua,{\it Phys.Lett.}
{\bf 328B}, 55 (1994).

\bibitem{7.} W.Pusz and S.L.Woronowicz,{\it Rep.Math.Phys.}{\bf 27}, 231
(1989);
  D.Fairlie and C.Zachos ,
{\it Phys.Lett.}{\bf 256B},43 (1991); D.Fairlie and J.Nuyts,
 {\it Z.Phys.C }
  {\bf 56}, 237 (1992);
 R. Chakrabarti and R. Jagannathan,
{\it J. Phys.A :Math.Gen.} {\bf 24}, L711 (1991); R.Jagannathan et. al.,
{\it J. Phys.A :Math.Gen.} {\bf 25}, 6429 (1992) and references therein.

\bibitem{8} P.P.Kulish, {\it Phys.Lett.}{\bf 161A}, 50 (1990); J.Van der Jeugt,
{\it J. Phys.A :Math.Gen.}{\bf 26},2405 (1993); A.Kempf,
{\it J.Math. Phys.}{\bf 34}, 969 (1993).

\bibitem{9} S.Meljanac,M.Milekovi\c and A.Perica,
{\it  }{\bf 28}, 79 (1994); M.Dore\s i\c , S.Meljanac and
M.Milekovi\c ,
{\it Fizika}{\bf 3}, 57 (1994).

\bibitem{10} H.S.Green, {\it Phys.Rev.}{\bf 90}, 170  (1953);
O.W.Greenberg and
A.M.L.Messiah, {\it Phys.Rev.}{\bf 138B}, 1155 (1965); {\it J.Math.Phys.}{\bf
6},
500 (1965).

\bibitem{11.} For a review,see:Y.Ohnuki  and S.Kamefuchi,
 {\it Quantum Field Theory and Parastatistics}
( University of Tokio Press, Tokio, Springer, Berlin, 1982).

\bibitem{12.} A.Bhattacharyya et. a l., {\it Phys.Lett.} {\bf 224B}, 384
(1989).

\bibitem{13.} A.B.Govorkov, {\it Nucl.Phys.}{\bf 365B},381 (1991);
{\it Phys.Elem.Part.Atomic Nucl. }{\bf 4},1341 (1993)(in russian).

\bibitem{14.} O.W.Greenberg, {\it  Phys. Rev. Lett.}{\bf 64}, 705 (1990);
{\it Phys. Rev. D}{\bf 43}, 4111 (1991); G.S.Agrawal and S.Chaturvedi,
{\it Mod.Phys.Lett.}{\bf A7},2407 (1992);
R.Speicher, {\it Lett.Math.Phys.}{\bf 27}, 97 (1993); M.Bozejko and
R.Speicher,
{\it Math.Ann.}{\bf 300}, 97 (1994);
 {\it Commun. Math. Phys.} {\bf 137}, 519  ;
 P.E.T.Jorgensen and R.F.Werner, {\it Commun. Math. Phys.} {\bf 164}, 455
(1994).

\bibitem{15.} P.E.T.Jorgensen,L.M.Schmitt and R.F.Werner,
{\it Pacific J.Math.}
{\bf 165},131 (1994).

\bibitem{16.} S.Meljanac and A.Perica, {\it J. Phys.A :Math.Gen. }{\bf 27},
 4737 (1994);\\
 {\it Mod.Phys.Letters}{\bf A9}, 3293 (1994);
V.Bardek,S.Meljanac and A.Perica, {\it Phys.Lett.}{\bf 338B}, 20 (1994) .

\bibitem{17.}  J.M.Leinaas and J.Myrheim,{\it Nuovo Cim.}{\bf 37}, 1 (1977);
 F.Wilczek,
{\it  Phys. Rev. Lett.} {\bf 48}, 1144 (1982);
J.Myrheim,{\it Anyons} preprint University of Trondheim 1993 (unpublished);
V.Bardek,M.Dore\s i\c $\,$ and S.Meljanac, {\it  Phys. Rev.}{\bf D49}, 3059
(1994);
V.Bardek,M.Dore\s i\c $\,$ and S.Meljanac, {\it Int.J.Mod.Phys.  }{\bf A9},
4185
(1994);
 A.K.Mishra and G.Rajasekaran, {\it Mod.Phys.Lett.}{\bf A9}, 419 (1994).

\bibitem{18.} F.D.M.Haldane,{\it  Phys. Rev. Lett.} {\bf 67}, 937 (1991);
 D.Karabali and
V.P.Nair, {\it Many-Body States and Operator Algebra for Exclusion Statistics}
 preprint IASSNS-HEP-94/88
 (to appear in Nucl.Phys.B) and references therein.

\bibitem{19.} A.P.Polychronakos, {\it Phys. Rev. Lett.}{\bf 69}, 703 (1991);
L.Brink et.al., {\it Nucl.Phys.}{\bf 401B}, 591 (1993).

\end{thebibliography}
\end{document}